\documentclass{amsart}
\usepackage[cp1251]{inputenc}
\usepackage[english]{babel}
\usepackage{amsmath}
\usepackage{amssymb}
\usepackage{amsfonts}

\begin{document}

\thispagestyle{empty}

 \title[The solution of
fractal diffusion retrospective problems]{Hermite functions with discontinuous coefficients for the solution of
fractal diffusion retrospective problems
 }%

\author{Yaremko O.E.}

\address{Oleg Emanuilovich Yaremko,
\newline\hphantom{iii}Penza State University,
\newline\hphantom{iii}str. Lermontov, 37, 
\newline\hphantom{iii} 440038, Penza, Russia}
\email{yaremki@mail.ru}

\maketitle {\small

\footnote{ Yaremko
Oleg, E-mail address: yaremki@mail.ru\par }
\begin{quote}
\noindent{\sc Abstract. }
In this article we study the
retrospective inverse problem. The retrospective inverse problem consists of
in the reconstruction of a priori unknown initial condition of the dynamic
system from its known final condition. Existence and uniqueness of the
solution is proved.

\medskip


 \end{quote}

\section{Introduction}
\label{sec:introduction}

In this article we study the retrospective inverse problem. The
retrospective Inverse problem consists of in the reconstruction of a priori
unknown initial condition of the dynamic system from its known final
condition. The direct problem of heat conductivity is well-posed; the
inverse problem is not well-posed. In mathematics the vast majority of
inverse problems set not well-posed - small perturbations of the initial
data (observations) can correspond to arbitrarily large perturbations of the
solution. The French mathematician Jacques Hadamard in 1939 defined, the
problem is called correct or well-posed problem if a solution exists, the
solution is unique, the solution's behavior hardly changes when there's a
slight change in the initial condition. If at least one of these three
conditions is not fulfilled, problems are termed ill-posed or not
well-posed. The most often in the case of ill-posed problems of the third
condition are violated the condition of the stability of solutions. In this
case, there is a paradoxical situation: the problem is mathematically
generated, but the solution cannot be obtained by conventional methods. A
classic example of ill-posed problem is retrospective problem for heat
equation on the real axis. Mathematically retrospective problem leads to a
Fredholm integral equation of the first kind:
\begin{equation}
\label{eq1}
\int\limits_{-\infty }^\infty {\frac{1}{2\sqrt {\pi \tau } }e^
{-\frac{(x-\xi )^2}{4\tau }} \hat {f}(\xi )d\xi =\hat {u}(\tau ,x),}
\end{equation}
in which $\hat {f}(x)$- is the initial distribution of the temperature
field, $\hat {u}(t,x)$- is the distribution of the fields in the moment of
time t. As shown in [1], the solution of equation (\ref{eq1}) expressed by the
formula:
\begin{equation}
\label{eq2}
\hat {f}(x)=\frac{1}{\sqrt \pi }\sum\limits_{j=0}^\infty {\frac{\hat
{u}^{(j)}(0)}{(2\sqrt \tau )^{n+1}j!}H_j \left( {\frac{x}{2\sqrt \tau }}
\right)} .
\end{equation}
\section{Problem Statement}
In the inverse problem of heat conductivity the initial distribution of
sources is unknown. The initial distribution of sources generates the
specified temperature distribution in an infinite piecewise-homogeneous rod
In. Mathematical statement of the problem consists in finding a solution
separatist system (n+1) equations of parabolic type
\begin{equation}
\label{eq3}
\left( {\frac{\partial }{\partial t}-a_j^2 \frac{\partial ^2}{\partial x^2}}
\right)u_j (t,x)=0,t>0,x\in I_n ,
\end{equation}
on the initial conditions
\begin{equation}
\label{eq4}
\left. {u_j (t,x)} \right|_{t=0} =f_j (x),x\in I_n ,
\end{equation}
on the boundary conditions
\begin{equation}
\label{eq5}
\left. {u_1 } \right|_{x=-\infty } =0,\left. {u_{n+1} } \right|_{x=\infty }
=0
\end{equation}
and the coupling conditions
\begin{equation}
\label{eq6}
\left[ {\alpha _{m1}^k \frac{\partial }{\partial x}+\beta _{m1}^k }
\right]u_k =\left[ {\alpha _{m2}^k \frac{\partial }{\partial x}+\beta
_{m2}^k } \right]u_{k+1,}
\end{equation}
\[
x=l_k ,k=1,...,n;m=1,2,
\]
here $u(t,x)$ - unknown function,
\[
u(\tau ,x)=\sum\limits_{k=2}^n {\theta (x-l_{k-1} )\theta (l_k -x)u_k (\tau
,x)+}
\]
\[
+\theta (l_1 -x)u_1 (\tau ,x)+\theta (x-l_n )u_{n+1} (\tau ,x),
\]
$f(x)$-unknown initial condition,
\[
f(x)=\sum\limits_{k=2}^n {\theta (x-l_{k-1} )\theta (l_k -x)f_k (x)+}
\]
\[
+\theta (l_1 -x)f_1 (x)+\theta (x-l_n )f_{n+1} (x)
\]
 $\theta(x)$ - set function, $\alpha _{mi}^k ,\beta _{mi}^k ,\gamma _{mi}^k ,\delta _{mi}^k ,\quad -$ given
real number, in which the condition of unlimited solvability of the problem
considered fulfilled [2].

The solution to problem (\ref{eq3})-(\ref{eq6}) is of the form:
\begin{equation}
\label{eq7}
u_k (t,x)=\sum\limits_{s=1}^{n+1} {\int\limits_{l_s }^{l_{s+1} } {H_{ks}
(t,x,\xi )f_s (\xi )d\xi ,} }
\end{equation}
where $H_{ks} (t,x,\xi )k,s=1,...,n+1$- influence function [2] of the mixed
boundary value problem.

Retrospective problem for the heat equation in the case of infinite
piecewise-homogeneous rod consists in the determination of the unknown
initial distribution of sources $f(x)$, which generates the specified
temperature distribution $u(\tau ,x)$ in the moment of time $t=\tau $.

\section{Operators transform }
Method of operators transform  is used to solve the problem [2].
Necessary definitions from [2], [11], [12]. The direct $J:\hat {f}\to f$ f
and inverse $J^{-1}:f\to \hat {f}$ operators transform  are set
equalities:
\[
f(x)=\int\limits_{-\infty }^\infty {\varphi (x,\lambda )\left(
{\int\limits_{-\infty }^\infty {e^{-i\lambda \xi }\hat {f}(\xi )d\xi } }
\right)d\lambda ,}
\]
\[
\hat {f}(x)=\int\limits_{-\infty }^\infty {e^{-i\lambda \xi }\left(
{\int\limits_{-\infty }^\infty {\varphi ^\ast (\xi ,\lambda )f(\xi )d\xi } }
\right)d\lambda .}
\]
Here $\varphi (x,\lambda ),\varphi ^\ast (x,\lambda )$-are the
eigenfunctions [13], [14] of the direct and coupling Sturm--Liouville
problems for the Fourier operator in piecewise-homogeneous axis In.
Eigenfunction
\[
\varphi (x,\lambda )=\sum\limits_{k=2}^n {\theta (x-l_{k-1} )} \theta (l_k
-x)\varphi _k (x,\lambda )+
\]
\[
+\theta (l_1 -x)\varphi _1 (x,\lambda )+\theta (x-l_n )\varphi _{n+1}
(x,\lambda )
\]
is a solution of the system of separate differential equations
\[
\left( {a_m^2 \frac{d^2}{dx^2}+\lambda ^2} \right)\varphi _m (x,\lambda
)=0,x\in \left( {l_m ,l_{m+1} } \right);m=1,...,n+1,
\]
on the coupling conditions
\[
\left[ {\alpha _{m1}^k \frac{d}{dx}+\beta _{m1}^k } \right]\varphi _k
=\left[ {\alpha _{m2}^k \frac{d}{dx}+\beta _{m2}^k } \right]\varphi _{k+1}
,
\]
on the boundary conditions
\[
\left. {\varphi _1 } \right|_{x=-\infty } =0,\left. {\varphi _{n+1} }
\right|_{x=\infty } =0
\]
Similarly eigenfunction
\[
\begin{array}{l}
 \varphi ^\ast (\xi ,\lambda )=\sum\limits_{k=2}^n {\theta (\xi -l_{k-1}
)\theta (l_k -\theta )} \varphi ^\ast (\xi ,\lambda )+ \\
 \theta (l_1 -\xi )\varphi _1^\ast (\xi ,\lambda )+\theta (\xi -l_n )\varphi
_{n+1}^\ast (\xi ,\lambda ) \\
 \end{array}
\]
is a solution of the system of separate differential equations
\[
\left( {a_m^2 \frac{d^2}{dx^2}+\lambda ^2} \right)\varphi _m^\ast (x,\lambda
)=0,x\in (l_m ,l_{m+1} );m=1,...,n+1,
\]
with the coupling conditions
\[
\frac{1}{\Delta _{1,k} }\left[ {\alpha _{m1}^k \frac{d}{dx}+\beta _{m1}^k }
\right]\varphi _k^\ast =\frac{1}{\Delta _{2,k} }\left[ {\alpha _{m2}^k
\frac{d}{dx}+\beta _{m2}^k } \right]\varphi _{k+1}^\ast ,x=l_k ,
\]
where
\[
\Delta _{i,k} =\det \left( {\begin{array}{l}
 \alpha _{1i}^k \beta _{1i}^k \\
 \alpha _{2i}^k \beta _{2i}^k \\
 \end{array}} \right)k=1,...,n;i,m=1,2,
\]
on the boundary conditions
\[
\left. {\varphi _1 } \right|_{x=-\infty } =0,\left. {\varphi _{n+1} }
\right|_{x=\infty } =0
\]
Let for some $\lambda $ of the considered boundary value problems have
nontrivial solutions $\varphi (x,\lambda ),\varphi ^\ast (x,\lambda )$, in
this case the number $\lambda $ is called the eigenvalue [13], [14],
corresponding solutions $\varphi (x,\lambda ),\varphi ^\ast (x,\lambda )$ -
is called the eigenfunctions of the direct and coupling Sturm--Liouville
problems, respectively. In the further we shall adhere to the following
normalization of eigenfunctions:
\[
\varphi _{n+1} (x,\lambda )=e^{ia_{n+1}^{-1} x\lambda },\varphi _{n+1}^\ast
(x,\lambda )=e^{-ia_{n+1}^{-1} x\lambda }.
\]
\section{Analogues of the system the Hermite functions on piecewise-homogeneous real
axis}
Let define analogues of the system the Hermite functions on
piecewise-homogeneous real axis:
\[
H_{j,n} (x)=\int\limits_{-\infty }^\infty {\varphi (x,\lambda )H_j \left(
{\frac{\lambda }{2\sqrt \tau }} \right)d\lambda ,}
\]
\[
H_{_{j,n} }^\ast (x)=\int\limits_{-\infty }^\infty {\varphi ^\ast (x,\lambda
)H_j \left( {2\sqrt \tau \lambda } \right)d\lambda .}
\]
where $H_j $- the system of the classical orthogonal Hermite functions [1].

\textbf{Lemma 1.} Functions $H_{j,n} (x),H_{_{j,n} }^\ast (x)$ form
biorthogonal system of functions by piecewise-homogeneous real axis.

\textbf{Proof.} We have the equality:
\[
\begin{array}{l}
 \int\limits_{-\infty }^\infty {H_{j,n} (x)} H_{_{k,n} }^\ast (x)dx= \\
 =\int\limits_{-\infty }^\infty {\left( {\int\limits_{-\infty }^\infty
{\varphi (x,\lambda )H_j (\lambda )d\lambda } } \right)} \left(
{\int\limits_{-\infty }^\infty {\varphi ^\ast (x,\beta )H_k (\beta )d\beta }
} \right)dx. \\
 \end{array}
\]
We change the integrals of places, we get:
\[
\begin{array}{l}
 \int\limits_{-\infty }^\infty {H_{j,n} (x)} H_{_{k,n} }^\ast (x)dx= \\
 =\int\limits_{-\infty }^\infty {H_j (\lambda )\left( {\int\limits_{-\infty
}^\infty {\varphi (x,\lambda )} \left( {\int\limits_{-\infty }^\infty
{\varphi ^\ast (x,\beta )H_k (\beta )d\beta } } \right)dx} \right)} d\lambda
. \\
 \end{array}
\]
On the decomposition theorem, we have:
\[
H_k (\lambda )=\int\limits_{-\infty }^\infty {\varphi (x,\lambda )} \left(
{\int\limits_{-\infty }^\infty {\varphi ^\ast (x,\beta )H_k (\beta )d\beta }
} \right)dx.
\]
Consequently,
\[
\int\limits_{-\infty }^\infty {H_{j,n} (x)} H_{_{k,n} }^\ast
(x)dx=\int\limits_{-\infty }^\infty {H_j (\lambda )H_k (\lambda )d\lambda }
=\delta _{j,k} .
\]
\section{Main result}
The problem of determining the initial distribution of the temperature field
$f(x)$ mathematically leads to the separate system of integral equations:
\begin{equation}
\label{eq8}
\sum\limits_{s=1}^{n+1} {\int\limits_{l_s }^{l_{s+1} } {H_{ks} (\tau ,x,\xi
)f_s (\xi )d\xi =u_k (\tau ,x)} };  k=1,...,n+1.
\end{equation}

Method of transformation operators applicable to solving separate system of
integral equations (\ref{eq8}).

\textbf{Theorem 1.} If the function $u(\tau ,x)\in {S}'(R)$ and for her the
condition
\[
e^{\tau ^2\lambda }(1+\lambda ^2)^{\alpha /2}\tilde {u}(\tau ,\lambda )\in
L_2 (R),
\]
that the separate system of integral equations (\ref{eq8}) has a unique solution
$f(x)\in H_2^\alpha (I_n )$ (definition $H_2^\alpha (I_n )$[6]), is
according to the formula:
\[
f(x)=\frac{1}{\sqrt \pi }\sum\limits_{j=0}^\infty {\frac{D_n
(u)}{2^jj!}H_{j,n} (x),}
\]
where
\[
D_n (u)=\frac{1}{2\pi }\int\limits_{-\infty }^\infty {(i\lambda )^j\tilde
{u}(\tau ,\lambda )d\lambda .}
\]
\textbf{Proof .} Let's apply the operators transform $J^{-1}$ to
separate system of integral equations (\ref{eq8}). As a result come to a model
integral equation (\ref{eq1}). Let's apply the operator $J$ in both parts of the
obtained equality (9); as a result, taking into account the continuity of
the operator$J$, we find the unknown distribution of temperature:
\[
f(x)=\frac{1}{\sqrt \pi }\sum\limits_{j=0}^\infty {\frac{\hat
{u}^{(j)}(0)\tau ^{\frac{j}{2}}}{\left( {2\sqrt \tau }
\right)^{n+1}j!}H_{j,n} (x).}
\]
Let's calculate numbers $\hat {u}^{(j)}(0).$
\[
\hat {u}^{(j)}(0)=\frac{1}{2\pi }\int\limits_{-\infty }^\infty {(i\lambda
)^j\left( {\int\limits_{-\infty }^\infty {e^{-i\lambda \xi }\hat {u}(\xi
)d\xi } } \right)d\lambda ,}
\]
from the definition of the operator conversion of $J$ the equality follows:
\[
\int\limits_{-\infty }^\infty {e^{-i\lambda \xi }\hat {u}(\xi )d\xi =}
\int\limits_{-\infty }^\infty {\varphi ^\ast (\xi ,\lambda )u(\xi )d\xi .}
\]
Thus,
\[
\hat {u}^{(j)}(0)=\frac{1}{2\pi }\int\limits_{-\infty }^\infty
{(i\lambda)}^{j}\tilde {u}(\tau ,\lambda )d\lambda .
\]
\section{Power function with discontinuous coefficients and its application}
We consider the Fourier transform of the Delta function [15]
\[
F\left[ {\delta (x)} \right]=\int\limits_{-\infty }^\infty {e^{-i\lambda
x}\delta (x)dx=1,}
\]
consequently
\[
F^{-1}\left[ 1 \right]=\frac{1}{2\pi }\int\limits_{-\infty }^\infty
{e^{-i\lambda x}dx=\delta (\lambda ).}
\]
We find as a consequence
\[
F\left[ {x^k} \right]=i^k\frac{d^k}{d\lambda ^k}\left( {\int\limits_{-\infty
}^\infty {e^{-i\lambda x}dx} } \right)=2\pi i^k\delta ^k(\lambda ).
\]
Let's define analog of the power function as follows
\[
x_n^k =(i)^k\left. {\frac{\partial ^k\varphi (x,\lambda )}{\partial \lambda
^k}} \right|_{\lambda =0}.
\]
We find equality from the definition of the operators transform
\[
J\left[ {x^k} \right]=x_n^k ;J^{-1}\left[ {x_n^k } \right]=x^k.
\]
The last equality means that the power function with discontinuous
coefficients is obtained by the action of the operators transform to the
power function.

\textbf{Theorem 2.} The ratio connects the generalized power function and
differentiation
\[
\frac{d^2}{dx^2}x_n^k =k(k-1)x_n^{k-2} .
\]
\textbf{Proof.} We have a chain of equalities
\[
\frac{d^2}{dx^2}\left( {x_n^k } \right)=\frac{d^2}{dx^2}J\left( {x^k}
\right)=J\left( {\frac{d^2}{dx^2}x^k} \right)=
J(k(k-1)x^{k-2} )=k(k-1)x_{n}^{k-2} .
\]
\section{Retrospective problem for the system of the diffusion equations}
Let's return to the solution of the separate system of integral equations
(\ref{eq8}) in space of the generalized functions [16] ${S}'$
\begin{equation}
\label{eq9}
f(x)=\frac{1}{2\pi }\int\limits_{-\infty }^\infty {e^{\lambda ^2\tau
}\varphi (x,\lambda )} \int\limits_\infty ^\infty {\varphi ^\ast (x,\xi
)u(\tau ,\xi )d\xi d\lambda .}
\end{equation}
We will convert the found solution. An analogue of the Taylor series for the
function $u(\tau ,\xi )$ is of the form
\begin{equation}
\label{eq10}
u(\tau ,\xi )=\sum\limits_{j=0}^\infty {\frac{u_j (\tau )}{j!}x_n^j .}
\end{equation}
From decomposition theorem have

\begin{equation}
\label{eq11}
u(\tau ,x)=\frac{1}{2\pi }\sum\limits_{s=1}^{n+1} {\int\limits_{-\infty }^\infty
{\varphi (x,\lambda )} \int\limits_{l_s }^{l_{s+1} } {\varphi _s^\ast (x,\xi
)u_s (\tau ,\xi )d\xi d\lambda .} }
\end{equation}
Let's find the decomposition of the generalized Taylor series for
eigenfunction from definition of the generalized power function
\[
\varphi (x,\lambda )=\sum\limits_{j=0}^\infty {\frac{(i\lambda )^j}{j!}x_n^j
.}
\]
Let's substitute this decomposition in a formula (\ref{eq10}) and let's integrate
term by term. We come to the formula (\ref{eq11}) in which
\begin{equation}
\label{eq12}
u_j (\tau )=\frac{1}{2\pi }\sum\limits_{s=1}^{n+1} {\int\limits_{-\infty
}^\infty {(i\lambda )^j} \int\limits_{l_s }^{l_{s+1} } {\varphi _s^\ast
(x,\xi )u_s (\tau ,\xi )d\xi d\lambda .} }
\end{equation}
Let's substitute decomposition of eigenfunction $\varphi (x,\lambda )$ in
the generalized power series in a formula (\ref{eq9}) we will receive
\[
f(x)=\frac{1}{2\pi }\int\limits_{-\infty }^\infty {e^{\lambda ^2\tau
}\varphi (x,\lambda )} \int\limits_\infty ^\infty {\varphi ^\ast (x,\xi
)\sum\limits_{j=0}^\infty {\frac{u_j (\tau )}{j!}\xi _n^j d\xi d\lambda =} }
\]
\[=\sum\limits_{j=0}^\infty {\frac{u_j (\tau )}{j!}} \frac{1}{2\pi
}\int\limits_{-\infty }^\infty {e^{\lambda ^2\tau }\varphi (x,\lambda )}
\int\limits_\infty ^\infty {\varphi ^\ast (\xi ,\lambda )\xi _n^j d\xi
d\lambda =}
\]
\[
=\sum\limits_{j=0}^\infty {\frac{u_j (\tau )}{j!}} \frac{1}{2\pi
}\int\limits_{-\infty }^\infty {e^{\lambda ^2\tau }\varphi (x,\lambda )}
2\pi i^j\delta ^j(\lambda )d\lambda =
\]
\[
=\sum\limits_{j=0}^\infty {\frac{u_j (\tau )(-1)^j}{j!}}
\int\limits_{-\infty }^\infty {\frac{\partial ^j}{\partial \lambda ^j}\left(
{e^{\lambda ^2\tau }\varphi (x,\lambda )} \right)} i^j\delta ^j(\lambda
)d\lambda =
\]
\[
=\sum\limits_{j=0}^\infty {\frac{u_j (\tau )}{j!}} H_{jn} (x),
\]
where designation is accepted
\[
\left. {H_{jn} (x)=(-i)^j\frac{\partial ^j}{\partial \lambda ^j}\left(
{e^{\lambda ^2\tau }\varphi (x,\lambda )} \right)} \right|_{\lambda =0} .
\]
Thus, the solution of the retrospective problem (\ref{eq9}) is obtained.

\textbf{Remark 1.} The generating function for the $H_{jn} (x)$ is the form
\[
e^{-\lambda ^2\tau }\varphi (x,-i\lambda )=\sum\limits_{j=0}^\infty
{\frac{H_{jn} (x)}{j!}\lambda ^j.}
\]
\textbf{Corollary 1}. If to choose $t=\frac{1}{2},n=1,\varphi (x,\lambda
)=e^{i\lambda x},$ then we obtain
\[
e^{-\frac{1}{2}\lambda ^2}e^{\lambda x}=\sum\limits_{j=0}^\infty
{\frac{H_{j1} (x)}{j!}\lambda ^j,}
\]
functions$H_{j1} (x)$-are the classical Hermite polynomials.

\section{Retrospective problem for fractal systems of diffusion equations}
Retrospective problem for fractal system of diffusion equations in the space
of generalized functions ${S}'$ leads to the separatist system of integral
equations:
\[
f(x)=\frac{1}{2\pi }\int\limits_{-\infty }^\infty {E_{\alpha ,1} (\lambda
^2\tau ^\alpha )\varphi (x,\lambda )} \int\limits_\infty ^\infty {\varphi
^\ast (x,\xi )u(\tau ,\xi )d\xi d\lambda ,}
\quad
k=1,...,n+1,
\]
where $E_{\alpha ,1} (z)$- the Mittag-Leffler function [9].

We get the solution of the fractal retrospective problem
\[
f(x)=\sum\limits_{j=0}^\infty {\frac{u_j (\tau )}{j!}H_{jn} (x),}
\]
where designation is accepted
\[
H_{jn} (x)\left. {=(-i)^j\frac{\partial ^j}{\partial \lambda ^j}\left(
{E_{\alpha ,1} \left( {\lambda ^2\tau ^\alpha } \right)\varphi (x,\lambda )}
\right)} \right|_{\lambda =0} .
\]
You can find the generating functions for $H_{jn} (x)$:
\begin{equation}
\label{eq13}
E_{\alpha ,1} \left( {-\lambda ^2\tau ^\alpha } \right)\varphi (x,-i\lambda
)=\sum\limits_{j=0}^\infty {\frac{H_{jn} (x)}{j!}\lambda ^j.}
\end{equation}
We find the explicit expression for the functions $H_{jn} (x)$. We find
decomposition of the left part of the formula (\ref{eq13}) in the Taylor series
\[
E_{\alpha ,1} \left( {-\lambda ^2\tau ^\alpha } \right)\varphi (x,-i\lambda
)=\sum\limits_{k=0}^\infty {\left( {-1} \right)^k\frac{\lambda ^{2k}\tau
^{\alpha k}}{\Gamma (k\alpha +1)}} \sum\limits_{m=0}^\infty {\frac{\lambda
^m}{m!}x_n^m =}
\]
\[
=\sum\limits_{j=0}^\infty {\lambda ^j} \sum\limits_{2k+m=j} {\left( {-1}
\right)^k\frac{\tau ^{\alpha k}}{\Gamma (k\alpha +1)m!}} x_n^m =
\]
\[
=\sum\limits_{j=0}^\infty {\lambda ^j} \sum\limits_{k=0}^{\left[
{\frac{n}{2}} \right]} {\frac{\left( {-1} \right)^k\tau ^{\alpha
k}j!}{\Gamma (k\alpha +1)(j-2k)!}} x_n^{j-2k} .
\]
We get the expression, comparing the two views $H_{jn} (x)$
\[
H_{jn} (x)=\sum\limits_{k=0}^{\left[ {\frac{n}{2}} \right]} {\frac{\left(
{-1} \right)^k\tau ^{\alpha k}j!}{\Gamma (k\alpha +1)(j-2k)!}} x_n^{j-2k} .
\]
If $n=1$, the formula takes the form
\[
H_{j1} (x)=\sum\limits_{k=0}^{\left[ {\frac{j}{2}} \right]} {\frac{\left(
{-1} \right)^k\tau ^{\alpha k}j!}{\Gamma (k\alpha +1)(j-2k)!}} x_{n}^{j-2k} .
\]
Define a fractal generalization of the Hermite polynomials
\[
H_j (x)=\sum\limits_{k=0}^{\left[ {\frac{j}{2}} \right]} {\frac{\left( {-1}
\right)^kj!}{\Gamma (k\alpha +1)(j-2k)!}} x_{n}^{j-2k} ,
\]
then the solution of the fractal retrospective problem has the form
\[
f(x)=\sum\limits_{j=0}^\infty {\frac{u_j (\tau )\tau ^{\beta j}}{j!}H_j
\left( {\frac{x}{\tau ^\beta }} \right),\beta =\frac{\alpha }{2}.}
\]
\textbf{Corollary 2}. In the hyperbolic case $\alpha =2,n=1$ the solution of
the retrospective problem has the form
\[
f(x)=\sum\limits_{j=0}^\infty {\frac{u_j (\tau )\tau ^j}{j!}H_j \left(
{\frac{x}{\tau }} \right),\beta =\frac{\alpha }{2}.}
\]
Where
\begin{equation}
\label{eq14}
H_j (x)=\sum\limits_{k=0}^{\left[ {\frac{j}{2}} \right]}
{\frac{(-1)^kj!}{((2k)!)(j-2k)!}x^{j-2k}}
=\frac{(1+x)^j+(1-x)^j}{2}.
\end{equation}
\textbf{Proof}. We replace $\lambda $ for $i\lambda $ in formula (\ref{eq14}) for
the solution of the direct problem. As a result, we get the formula
\[
u(\tau ,x)=\sum\limits_{j=0}^\infty {\frac{f_j \tau ^{\beta j}}{j!}H_j^\ast
\left( {\frac{x}{\tau ^\beta }} \right)} ,H_j^\ast (z)=i^jH_j (-iz).
\]
We use the formula (\ref{eq13})
\[
H_{j1} (x)=(-i)^j\frac{\partial ^j}{\partial \lambda ^j}\left. {(\cos
(\lambda \tau )e^{i\lambda x})} \right|_{\lambda =0} =
\]
\[
=\frac{1}{2}(-i)^j\frac{\partial ^j}{\partial \lambda ^j}\left.
{(e^{i\lambda (\tau +x)}+e^{i\lambda (\tau -x)})} \right|_{\lambda =0} =
\]
\[
=\frac{(\tau +x)^j+(\tau -x)^j}{2}.
\]
In the end we find the solution of the Cauchy problem for the hyperbolic
equation
\[
f(x)=\sum\limits_{j=0}^\infty {u^j(\tau ,0)} \frac{(\tau +x)^j+(\tau
-x)^j}{2j!}.
\]
Thus,
\[
f(x)=\frac{u(\tau ,\tau +x)+u(\tau ,\tau -x)}{2}.
\]
\textbf{Remark.} If, as an example, take
\[
u(t,x)=g(x+t-\tau )+g(x-t+\tau ),
\]
the $u(\tau ,x)=2g(x)$ and, so,
\[
f(x)=g(x+\tau )+g(x-\tau )=u(0,x).
\]
\section{The inverse d
Dirichlet problem for a half-plane}
Solution of the inverse Dirichlet problem for the right half-plane has the
form:
\begin{equation}
\label{eq15}
f(y)=Re\frac{1}{\pi }\int\limits_0^\infty {e^{\lambda l}e^{i\lambda y}}
\int\limits_{-\infty }^\infty {e^{-i\lambda \eta }u(l,\eta )d\eta d\lambda
.}
\end{equation}
Let's repeat the above reasoning. Let's receive expressions for analogues of
Hermite polynomials:
\[
H_{j1}(x)=Re(-i)^{j}\frac{\partial^{j}}{\partial\lambda^{j}}\left
(e^{\lambda l}e^{i\lambda y}) \right|_{\lambda =0} =
\]
\[
=Re(-i)^j(l+iy)^j=\frac{(y+li)^j+(y-li)^j}{2}.
\]
As a result for the solution of the inverse Dirichlet problem we receive the
representation in the form of the sum of the Taylor series:
\[
f(y)=\sum\limits_{j=0}^\infty
{\frac{u^j(l,0)}{j!}\frac{(y+li)^j+(y-li)^j}{2}} .
\]
\textbf{Corollary 3.} If the function $f(y)$ admit continued with the real
axis of the complex plane as a whole, the
\[
f(y)=Re u(l,y+li).
\]
\textbf{Example 2.} Let $u(x,y)=(x-l)^2-y^2$, then $u(l,y)=-y^2$ therefore,
we find
\[
f(y)=Re(-(y+li)^2)=l^2+y^2=u(0,y).
\]
\section{Conclusion}
In this article the formal solution of the retrospective problem is
provided. The third aspect in determining the well-posed problem is not
taken into account. Theorem of existence and uniqueness of solution are
given. From the analysis of the formula (9): the solution of the
retrospective heat problem with discontinuous coefficients is received by
replacement in the final result of the Hermite functions [1] on the Hermite
functions with discontinuous coefficients, defined in the article. The
derivatives $\hat {u}^j(0)$ are necessary to replace on $D_n (u)$. The
noticed analogy allows hoping on the possibility of obtaining the solution
of problems of mathematical physics in which the Hermite functions with
discontinuous coefficients.

\textbf{Yaremko Oleg Emanuilovich}, Physics and Mathematics Science (PhD),
associate professor, Managing chair of the mathematical analysis, Penza
state university, 440026, Russia, Penza, Lermontov's street, 37. E-mail
address: yaremki@mail.ru

\end{document}